\documentclass[twoside,12pt]{article}

\usepackage{amssymb,amsmath}

\newenvironment{proof}{\begin{trivlist}\item[]{\it
Proof.}}{\hfill$\square$\end{trivlist}}

\newtheorem{theorem}{Theorem}
\newenvironment{corollary}{\begin{trivlist}\item[]{\bf\noindent Corollary }\it}
{\end{trivlist}}

\newtheorem{lemma}[theorem]{Lemma}
\newtheorem{proposition}[theorem]{Proposition}

\newenvironment{remark}{\begin{trivlist}\item[]{\bf\noindent Remark }}
{\end{trivlist}}

\begin{document}

\date{}
\newcommand{\Z}{\mathbb Z}\newcommand{\C}{\mathbb C}\newcommand{\Q}{\mathbb Q}
\newcommand{\Pf}{{\mathit {Pf} }}
\newcommand{\adI}{^{(i)}}

\newcommand{\adIJ}{^{(ij)}}

\title{Mod 2 indecomposable orthogonal invariants}

\author{M. Domokos ${}^{a}$ 
\thanks{Corresponding author. 
Supported through a European Community Marie Curie Fellowship, 
held at the University of Edinburgh. 
Partially supported by OTKA grants No. T 34530 and T 046378.}\; 
and P. E. Frenkel ${}^{b}$ 
\thanks{Supported through a European Union Marie Curie Studentship held at the
University of Warwick. Partially supported by OTKA grants T 042769 and T 046365.}
\\ 
\\ 
{\small ${}^a$ R\'enyi Institute of Mathematics, Hungarian Academy of 
Sciences,} 
\\ {\small P.O. Box 127, 1364 Budapest, Hungary,} 
{\small E-mail: domokos@renyi.hu } 
\\   
\\ 
{\small ${}^b$ Institute of Mathematics, Budapest Technical University,} 
\\ {\small P.O.B. 91,  1521 Budapest, Hungary,}  
\\ {\small E-mail: frenkelp@renyi.hu}  
}
\maketitle 
\begin{abstract} 
Over an algebraically closed base field $k$ of characteristic 2,
the ring  $R^G$ of invariants is studied,  $G$ being the orthogonal group
 $O(n)$ or the special orthogonal group $SO(n)$ and acting naturally
 on the coordinate ring $R$ of the $m$--fold direct sum
 $k^n\oplus\cdots\oplus k^n$ of the standard vector representation. 
It is proved for $O(n)$ $(n\geq 2)$ and for $SO(n)$ $(n\geq 3)$ 
that there exist $m$--linear invariants with $m$ arbitrarily large 
that are indecomposable (i. e.,  not expressible as  polynomials in  
invariants of lower degree).  In fact, they are explicitly 
constructed for all possible values of $m$. Indecomposability of 
corresponding invariants over $\Z$ immediately   follows.  
The constructions rely on analysing the 
Pfaffian 
of the skew-symmetric matrix whose entries above the diagonal are the scalar products of  the vector variables. 
\end{abstract} 



\section{Introduction and notation}\label{Prelim}
\subsection{The orthogonal group}\label{The-orth}
Let $F$ stand for an algebraically closed field of arbitrary characteristic.
Denote coordinates in $F^n$ by $x_1$, $y_1$,  \dots, $x_\nu$,
$y_\nu$ if $n=2\nu$ or  by $x_1$, $y_1$, \dots, $x_\nu$, $y_\nu$,
$z$ if $n=2\nu+1$.
The {\emph {standard quadratic form}} $q:F^n\to F$ is
\[q\overset{\rm {def}}=x_1y_1 +\dots +x_\nu y_\nu
\mbox{ when }n=2\nu,\]
and
\[q\overset{\rm {def}}=x_1y_1 +\dots +x_\nu y_\nu+z^2
\mbox{ when }n=2\nu+1.\]
The \emph{orthogonal group} $O(n,F)$ is defined as the
group of linear isomorphisms of $F^n$ that leave the quadratic form
$q$ invariant.  The \emph{special orthogonal group}  $SO(n,F)$ is defined as
the component of $O(n,F)$ containing the identity.
Recall that the polar form $\beta$ of $q$ is the symmetric bilinear form
on $F^n$
given by
\[\beta\left(v^{(1)},v^{(2)}\right)\overset{\rm{def}}=
q\left(v^{(1)}+v^{(2)}\right)-q\left(v^{(1)}\right)-q\left(v^{(2)}\right).\]
The non-degeneracy of $q$ means that $\beta(v,\cdot)=0$  and $q(v)=0$
together imply $v=0$.
Note that up to base change, $q$ is the only non-degenerate quadratic form on
$F^n$.

Throughout this paper $k$ stands for an algebraically closed field of
characteristic 2, and we write just $O(n)$ and $SO(n)$ when $k$ is
to be understood. All elements of $O(n)$ have determinant $1$.
The algebraic group $O(n)$ is
connected for odd $n$ and has two components for even $n$. 
Thus, $SO(2\nu+1)=O(2\nu+1)$, whereas $SO(2\nu)$ is a subgroup of
index 2 in $O(2\nu)$.

\subsection{Invariants}
We write $R$ or $R_{n\times m}$ for the $F$--algebra of polynomials in
the coordinates of the indeterminate $n$--dimensional vectors
$v^{(1)}$, \dots, $v^{(m)}$.  A letter $G$ in the superscript indicates
the subalgebra  formed by the polynomials invariant under
the group $G=O(n,F)$ or $G=SO(n,F)$ acting on $m$--tuples of vectors in the
obvious way. 

Next we recall from \cite{df} some distinguished elements in $R$. Set
\[Q^{(i)}\overset{\rm{def}}=q\left(v^{(i)}\right)
\mbox{ and }
B^{(ij)}
\overset{\rm{def}}=\beta\left(v^{(i)},v^{(j)}\right)\]
for $1\leq i,j\leq n$.
More explicitly, for $n=2\nu$ we have
\[B^{(ij)}=
x_1^{(i)}y_1^{(j)}+y_1^{(i)}x_1^{(j)}+ \dots
+x_\nu^{(i)}y_\nu^{(j)}+y_\nu^{(i)}x_\nu^{(j)},\]
whereas for $n=2\nu+1$ we have
\[B^{(ij)}=
x_1^{(i)}y_1^{(j)}+y_1^{(i)}x_1^{(j)}+ \dots
+x_\nu^{(i)}y_\nu^{(j)}+y_\nu^{(i)}x_\nu^{(j)}
+2z^{(i)}z^{(j)}.\]
Let
\[D^{(i_1, \dots , i_n)}\overset{\rm{def}}=\det \left[v^{(i_1)},
\dots , v^{(i_n)}\right]\]
be the determinant of the matrix that
has $v^{(i_1)}$, \dots, $v^{(i_n)}$ as its columns.
 Then  $
Q^{(i)}$ and $B^{(ij)}$  are
orthogonal invariants,
and  $D^{(i_1, \dots , i_n)}$ is  a  
special orthogonal invariant.

The classical ``first fundamental theorem" for the (special) orthogonal
 group asserts that when $F$ is of characteristic zero,
 the algebra $R^{O(n,F)}$ is generated by the
 scalar products $B^{(ij)}$ of the
 indeterminate vectors under consideration, and the algebra
 $R^{SO(n,F)}$ is generated by the scalar products and the
 determinants.
 This has been discussed along
 with the analogous results for the other  classical groups by
 Hermann Weyl in \cite{W}. 
 De Concini and Procesi \cite{CP} gave a characteristic free treatment 
to the subject, in particular, they proved that the first fundamental 
theorem for the (special) orthogonal group remains unchanged 
in odd characteristic. 
Concerning characteristic $2$, Richman \cite{Ri} proved later 
that the algebra $R^G$ for the group $G$ 
preserving the bilinear form 
$x_1^{(1)}x_1^{(2)}+\cdots+x_n^{(1)}x_n^{(2)}$ 
is generated in degree 
$1$ and $2$. 
However, though this group preserves the quadratic form 
$x_1^2+\cdots+x_n^2$, it is not the so-called `orthogonal group' 
in characteristic $2$: the quadratic form $x_1^2+\cdots+x_n^2$ 
is the square of a linear form, hence is degenerate. So in characteristic $2$ the question
about vector invariants of the orthogonal group remains open,
and was addressed in \cite{df} (motivated by an observation from
\cite{DKZ}).
In particular, it was shown in \cite{df} that the field of rational
$O(n)$--invariants is generated by the obvious quadratic invariants
in characteristic $2$ as well.
However, the behaviour of polynomial invariants turned out to be very much
different, see Section~\ref{high}.

%%%%%%%%%%%%%%%%%%%%%%%%%%%%%%%%%%%%%%%%%%%%%%%%%%%%%%%%%%%%%%%%

\section{Constructing indecomposable invariants}\label{high}
An element of the polynomial algebra $R$ is called \emph{$m$--linear}  if it is multilinear in the vector 
variables $v^{(1)}$, \dots, $v^{(m)}$. 
An element of the  algebra $R^G$ of invariants is said to be an
\emph{indecomposable} 
element if it is not contained in the subalgebra of $R^G$ generated by the 
elements of lower degree.
When $F=\C$, we  write  $R(\Z)=R_{n\times m}(\Z)$ for the subring of $R$
consisting of polynomials with integer coefficients, and $R^G(\Z)$ for the ring
$R^G\cap R(\Z)$ of $G$--invariants with integer coefficients  (note that $G$
here is still the complex group $O(n,\C)$ or $SO(n,\C)$).
In this case we say that an element of $R^G(\Z)$ is an indecomposable element
if it is not contained in the subring of $R^G(\Z)$ generated by the elements
of lower degree.

Let  $n=2\nu$ or $n=2\nu+1$. Over an algebraically closed base field $k$ of characteristic 2, we
construct an  indecomposable $m$--linear $O(n)$--invariant if  
 $m>2\nu\geq 2$ 
%are both even  or if $m\geq n\geq 3$ with $n$ 
and an   indecomposable $m$--linear $SO(n)$--invariant if $m\geq n\geq 3$, except when $n$ is even and $m$ is 
odd (in which case there are no $SO(n)$--invariants of degree $m$ at all ---
either see \cite[Theorem 4.5(ii)]{df} or just consider the action of the $\nu$--dimensional torus $SO(2)^\nu$ whose elements are the diagonal matrices in $SO(2\nu)$). 
%   For $m>n\geq 2$ both even, we also construct 
%an  indecomposable $m$--linear $O(n)$--invariant. 
All shall be constructed as modulo 2 images of invariants with integer coefficients, which therefore 
have the  above indecomposability properties 
 over $\Z$.

This proves our conjecture formulated in \cite{df}. Constructions     for
$n\leq 4$ were given there. The paper \cite{DKZ} contained a
more sophisticated proof of the existence of high-degree indecomposable
invariants in  the $SO(4)$ case.

We shall use the symbol $*$ to mean any one of the two letters $x$
and $y$.
We define the sign ${\rm sgn}\; \rho$ of  
an $m$--linear monomial 
\begin{equation}\label{monom} 
\rho =z^{(i_{01})}\cdots z^{(i_{0m_0})}
*_1^{(i_{11})} \cdots *_1^{(i_{1m_1})}\cdots\cdots *_\nu^{(i_{\nu 1})}\cdots *_\nu^{(i_{\nu m_\nu})} 
\end{equation}
with $i_{r1}<\cdots <i_{rm_r}$ for each $r$ to be the
sign of the permutation
$i_{01}$, \dots, $i_{0m_0}$, $i_{11}$, \dots, $i_{1m_1}$, \dots\dots, $i_{\nu 1}$, \dots, $i_{\nu m_\nu}$ of the indices 1, \dots, $m$.

In our first two propositions, which form the technical heart of the paper, we shall be working over $\Z$.  We write
$$\Pf ^{(i_1,\dots,i_{2\mu})}={\mathrm {Pf}}\left(\begin{matrix} 
0 & B^{(i_1i_2)} & B^{(i_1i_3)} & \dots & B^{(i_1|i_{2\mu})} \\ 
& 0 & B^{(i_2i_3)} & \dots  & B^{(i_2|i_{2\mu})}\\  
&&0&\dots &B^{(i_3|i_{2\mu})}\\&&&\ddots&\vdots\\ &&&&0\end{matrix}\right)\in 
R^{O(n,\C)}(\Z)$$
 for the Pfaffian of the 
$2\mu\times 2\mu$ skew-symmetric matrix  
whose upper half consists of the $B$'s corresponding to the indices 
$i_1$, \dots, $i_{2\mu}$. It  is invariant under $O(n,\C)$.  
It is multilinear in the vectors $v^{(i_1)}$, \dots,  $v^{(i_{2\mu})}$ if there is no repetition in the upper indices.
(See e. g. \cite[Appendix B.2]{gw} for the definition and basic 
properties of Pfaffians.)

\begin{proposition}\label{Bszorzat}
The coefficient in $\Pf ^{(1,\dots,2\mu)}$    of a $2\mu$--linear
monomial $\rho$ as in (\ref{monom}) above is zero unless, for each $r$,   $m_r=2\mu_r$ 
is even and the  $x_r$ and the $y_r$ occur in an alternating order. 
In this case, the coefficient is $$2^{\mu -
|\{r>0:\mu_r\geq 1\}|
}
{\rm sgn}\;\rho.$$
\end{proposition}

\newcommand{\sgn}{{\rm{sgn}}\;}

\begin{proof}
The coefficient of $\rho$ is given by substituting the corresponding  vectors of the standard basis into $\Pf ^{(1,\dots,2\mu)}$.  This amounts to calculating 
the Pfaffian of the skew-symmetric matrix given by the $B$'s of
these vectors.  Permuting rows and columns turns it into a direct sum of matrices and shows that the Pfaffian  is $\sgn\rho$ times a product of $m_r\times m_r$ Pfaffians, one factor for each  index $r$. 
So, assuming that   $\rho$ has
a non-zero coefficient,  each $m_r$ must be  even.  Set $\mu_r=m_r/2$.

If, for some $r>0$, we have two occurrences of  $x_r$ with no  $y_r$
in between (or {\it vice versa}), then the $r$th 
Pfaffian has two identical (adjacent) rows, which makes it zero.   So the $x_r$ and $y_r$ must occur in an alternating order for each $r>0$.
Therefore, the $r$th Pfaffian, for $r>0$, is that of the $m_r\times m_r$ 
checkerboard matrix with entries $b^{(ij)}={\rm sgn} (j-i)$ for odd $j-i$ and zero otherwise.   
Expanding by the first row and using 
induction on $\mu_r$ shows that this is $2^{\mu_r-
1
}$ if $\mu_r\geq 1$.  It is $1=2^{\mu_r}$ if $\mu_r=0$.

For odd $n$ and  $r=0$, we have the Pfaffian of the $m_0\times m_0$ matrix
with entries $b^{(ij)}=2{\rm sgn}(j-i)$.  Expanding by the first row and using
 induction on $\mu_0$ shows that this is $2^{\mu_0}$.  

Adding up the exponents of 2 yields the result.
\end{proof}

\begin{proposition}\label{DBszorzat}
Let $
m\geq n$ with $m\equiv n\mod 2$, and set $\mu=\lfloor m/2\rfloor$.
Consider the $m$--linear polynomial
 \begin{equation}\label{DB}
\sum
\sgn \pi\cdot D^{(\pi (1),\dots, \pi (n))}\Pf ^{(\pi (n+1),\dots,\pi(m))},
\end{equation} 
the sum being extended over
those permutations $\pi\in\mathfrak S_m$   that satisfy 
$\pi (1)<\cdots <\pi (n)$ and  $\pi(n+1)<\cdots <\pi(m)$.
The coefficient in (\ref{DB}) of an $m$--linear
monomial $\rho$ as in (\ref{monom}) above is zero unless, for each $r>0$,   $m_r$ is even and strictly 
positive, and the  $x_r$ and the $y_r$ occur in an alternating order. In this case, the coefficient is 
$$2^{\mu-\nu}(-1)^{\left|\left\{r>0 \; : \;*_r^{(i_{r1})}=y_r^{(i_{r1})}\right\}\right|}\sgn\rho.$$ 
\end{proposition}

\begin{proof}
Assuming that $\rho$ has a non-zero coefficient, it follows trivially that for each 
$r>0$, $m_r$ must be even and strictly positive.  Set $\mu_r=m_r/2$.   
We call a choice of one $x_r$ and one $y_r$ from $\rho$ a {\emph {good}} choice if the remaining $2(\mu_r-1)$ of the $x_r$ and $y_r$ 
occur in an alternating order.  Any choice 
of one $z$ from $\rho$ shall be called a good choice. Now  Proposition~\ref{Bszorzat}  
shows that  the terms in the sum 
(\ref{DB}) in which $\rho$ has non-zero coefficients correspond to the simultaneous good choices (for each $r$) of one of each letter from $\rho$, and the coefficient of $\rho$ in such a term 
is $$2^{\mu-\nu-
|\{r>0:\mu_r\geq 2\}|
}\sgn\rho\cdot\prod_r \pm 1,$$
where the $r$th $\pm 1$ is the sign of that permutation of the letters 
$x_r$ and $y_r$ in $\rho$   (resp. the letters $z$ in $\rho$)
which puts the chosen one(s) in front and in alphabetical order, leaving 
the rest in the order they had in $\rho$.

It follows that the coefficient of $\rho$ in (\ref{DB}) is
$$2^{\mu-\nu-
|\{r>0 : \mu_r\geq 2\}|
}\sgn\rho\cdot\prod_r\sum \pm 1,$$
where the $r$th summation is  over the good  choices of one $x_r$ and one $y_r$  (resp. one $z$)  
from $\rho$.

If we have two occurrences of $x_r$ with no $y_r$ in between (or {\it vice 
versa}), then the $r$th sum is zero, for we can pair off the choices by interchanging the r\^ole of the two adjacent $x_r$'s, and the two $\pm 1$'s  in each pair will 
cancel.  So the $x_r$ and $y_r$ must occur in an alternating order for each $r>0$.  In this case, the $r$th sum is
$\pm 1$ if  $\mu_r=1$   and $\pm\left(1+\sum_{j=1}^{m_r-1}(-1)^{j-1}\right)=\pm 2$ if $\mu_r\geq 2$.
The sign is $+$ if $x_r$ comes first and $-$ if $y_r$ comes first.  

For odd $n$, the 0th sum  is $\sum_{j=1}^{m_0}(-1)^{j-1}=1$. 

Adding up the exponents of 2 and $(-1)$ respectively, we arrive at the result.
\end{proof}

For $\mu\geq\nu$, divide the Pfaffian  in Proposition~\ref{Bszorzat} 
by $2^{\mu-\nu}$ to get a $2\mu$--linear $O(n,\C)$--invariant with integer coefficients, call it $\tilde g=\tilde g
_{n\times 2\mu}\in R_{n\times 2\mu}^{O(n,\C)}(\Z)$.  
View $\tilde g$ modulo 2  to get a 
$2\mu$--linear $O(n)$--invariant with coefficients in the prime subfield 
$\mathbb F_2$ of $k$, 
call it $g=g_{n\times 2\mu}\in R_{n\times 2\mu}^{O(n)}$. 
Analogously, for $m\geq n$ with $m\equiv n\mod 2$, divide the sum (\ref{DB}) of Proposition~\ref{DBszorzat} by $2^{\mu-\nu}$ to get an $m$--linear $SO(n,\C)$--invariant   
with integer coefficients, call it $\tilde h=\tilde h_{n\times m}\in R_{n\times m}^{SO(n,\C)}(\Z)$. View $\tilde h$ modulo 2  to get an 
$m$--linear $SO(n)$--invariant with coefficients in $\mathbb F_2$, call it $h=h_{n\times m}\in R_{n\times m}^{SO(n)}$. 
Note that invariance over $k$ in both cases follows from that over $\C$ by \cite[Lemma 3.2]{df}.

By
Proposition~\ref{Bszorzat}, $g$ is the sum of those $m$--linear monomials 
$\rho$ that, when written in the form \eqref{monom}, have strictly positive and even $m_r$ 
for all $r>0$, and the $x_r$ and $y_r$ occur in an alternating order. By Proposition~\ref{DBszorzat}, the same holds for $h$.
This shows in particular that for $m\geq n$ both even, $g=h$.  It follows that  $h$ is an $O(n)$--invariant whenever defined --- remember that $SO(2\nu+1)=O(2\nu+1)$.  

\bigskip

For odd $n$ and any $m\geq n$, exactly one of $g$ and $h$ is defined. 
We shall  prove that it is an indecomposable $SO(n)$--invariant if $n\geq 3$.
We define the \emph{multiplicity in $x$, $y$} of a polynomial 
to be the minimum of the total degrees of its monomials in the $x$, $y$ 
variables. (For  a homogeneous polynomial, this is  the difference 
between the degree and the degree in the $z$ variables.)
We observe that $g$ and  $h$
 have the smallest possible multiplicity in $x$, $y$.
Indeed, their multiplicity in $x$, $y$ is $2\nu$, and we have

\begin{lemma}\label{codeg}
The multiplicity in $x$, $y$
of any $ m$--linear $SO(2\nu+1)$--invariant is 
at least $\min (m,2\nu)$.
\end{lemma}

\begin{proof}
For $m\leq 2\nu$ we have
\cite[Theorem 4.9]{df}  that says that  the algebra $R_{(2\nu+1)\times
m}^{SO(2\nu+1)} $ is generated by the $Q\adI$ and the $B\adIJ$, which do not
involve the $z$ variables, so the 
 multiplicity in $x$, $y$ of any $m$--linear invariant
is $m$.  The case $m\geq 2\nu$ is reduced to this as follows.

Indirectly assume that an $m$--linear $SO(2\nu+1)$--invariant has a monomial
of the form $*_{r_1}^{(1)}\cdots *_{r_d}^{(d)}z^{(d+1)}\cdots z^{(m)}$ with
$0\leq d<2\nu$. Recall that  the $z$ axis in the space $k^{2\nu+1}$ 
is the radical of the symmetric bilinear form $\beta$, so its unit vector $e$ is stabilised by  
$SO(2\nu+1)$.   It follows that substituting $e$ 
for the vector variables $v^{(d+2)}$,\dots, $v^{(m)}$
in our $SO(2\nu+1)$--invariant yields a $(d+1)$--linear
$SO(2\nu+1)$--invariant, having a monomial of the form  $*_{r_1}^{(1)}\cdots
*_{r_d}^{(d)}z^{(d+1)}$. But $d+1\leq 2\nu$, so we have  a contradiction.
\end{proof}

Since the 
multiplicity in $x$, $y$ 
of a product of homogeneous polynomials is the  sum of the 
multiplicities of the factors, 
it follows for $\nu\geq 1$ that the multiplicity in $x$, $y$ 
of any decomposable  
$m$--linear invariant with $m>2\nu$ is strictly
greater than $2\nu$, so we get

\begin{theorem}\label{odd}
Let $n\geq 3$ be odd and $m\geq n$.  Then the $m$--linear $SO(n)$--invariant 
$g$ or $h$ (whichever one is defined) 
is indecomposable. 
\end{theorem}

\begin{corollary}
Let $n\geq 3$ be odd and $m\geq n$.  
Then the $m$--linear invariant $\tilde g$ or $\tilde h$ (whichever one is defined) 
is 
indecomposable in the ring  $R^{SO(n,\C)}(\Z)$.
 When $m$ is even, $\tilde g$ is, {\rm a fortiori}, indecomposable also 
in the ring  $R^{O(n,\C)}(\Z)$. 
\end{corollary}

Another, independent proof of Theorem~\ref{odd} for $m>n$ 
will be given in a remark following the
discussion of the $O(2\nu)$ case. 
Note that the case $m=n$
is immediate from \cite[Theorem 4.9]{df} cited in the proof above.

\medskip

We now wish to
prove for $n=2\nu$ and $\mu>\nu$ that the $O(2\nu)$--invariant $g$ is
indecomposable.  We pass to the new linear coordinates $t_r=x_r$ and
$s_r=x_r+y_r$ $(r=1,\dots ,\nu)$.  We define
the \emph{multiplicity in $s$} of a  polynomial to be the minimum of the total
degrees of its monomials in the $s$ 
variables, and we observe that $g$ 
 has the smallest possible  multiplicity in $s$.
Indeed, we have

\newcommand{\adi}{^{(1)}}
\newcommand{\adJ}{^{(j)}}

\begin{lemma}\label{mult}
The multiplicity  in $s$ of any $ 2\mu$--linear $O(2\nu)$--invariant
$p=p\left(v\adi,\dots, v^{(2\mu)}\right) $ is 
at least $\min (\mu,\nu)$.
\end{lemma}

\begin{proof}
For $\mu\leq \nu$ we have
\cite[Theorem 4.9]{df}  that says that  the algebra $R_{2\nu\times
2\mu}^{O(2\nu)} $ is generated by the $Q\adI$, which are quadratic in the
corresponding $v\adI$, and by the
$$B\adIJ=\sum_{r=1}^\nu\left(t_r\adI\left(t_r\adJ+s_r\adJ\right)+\left(t_r\adI+s_r\adI\right)t_r\adJ\right)
=\sum_{r=1}^\nu\left(t_r\adI s_r\adJ +s_r\adI t_r\adJ\right),$$
which are bilinear in $v\adI$, $v\adJ$ and have multiplicity 1 in $s$, so the multiplicity  in $s$
of any $2\mu$--linear invariant $p$
is $\mu$.

Now let  $\mu\geq \nu$.  It suffices  to prove that if $0\leq d<\nu$, then the
coefficient in $p$ of the monomial
\begin{equation}\label{monomial} 
s_{r_1}^{(1)}\cdots s_{r_d}^{(d)}t_{r_{d+1}}^{(d+1)}\cdots
t_{r_{2\mu}}^{(2\mu)}
\end{equation}
 is zero. 
This coefficient, expressed using
the original $x$, $y$ coordinates, 
is the sum of the coefficients 
in $p$ of the $2^{2\mu-d}$ monomials
$y_{r_1}^{(1)}\cdots y_{r_d}^{(d)}*_{r_{d+1}}^{(d+1)}\cdots
*_{r_{2\mu}}^{(2\mu)}$.  If  $\{r_{d+1},\dots,r_{2\mu}\}\not\subseteq\{r_1,\dots,r_d\}$, say $r_{d+1}\not\in\{r_1,\dots,r_d\}$, then these monomials can be paired off via the reflection $x_{r_{d+1}}\leftrightarrow y_{r_{d+1}}$, and the two coefficients in each pair are equal due to the $O(2\nu)$--invariant 
property of $p$.  

Now suppose that   $\{r_{d+1},\dots,r_{2\mu}\}\subseteq\{r_1,\dots,r_d\}$.
Since $d<2\mu-d$, the number of occurrences of at least one of the indices
$1$,\dots, $\nu$
among $r_1$,\dots , $r_d$ is less than among $r_{d+1}$, \dots ,
$r_{2\mu}$.  We may assume  $$r_1=\dots=r_a=1\ne r_{a+1},\dots, r_d$$ and
$$r_{d+1}=\dots=r_{d+a+1}=1.$$  Consider the $O(2\nu)$--invariant
$$p\left(u^{(r_1)},\dots,u^{(r_d)}, v^{(d+1)},\dots, v^{(d+a+1)},
w^{(r_{d+a+2})}, \dots, w^{(r_{2\mu})}\right)$$
depending on a new set
of vector variables whose cardinality is
$$|\{r_1,\dots,r_d\}|+a+1+|\{r_{d+a+2},\dots,r_{2\mu}\}|
\leq d-a+1+a+1+d=2(d+1)\leq 2\nu.$$
By \cite[Theorem 4.9]{df}
again, this invariant must be  a polynomial in the $q$'s
and $\beta$'s of its vector variables.  As it is linear in each of
$v^{(d+1)}$,\dots, $v^{(d+a+1)}$, these can 
be involved only via their $\beta$'s with each other or with 
the $u$'s and $w$'s.  
To get
the coefficient of the monomial \eqref{monomial}, we substitute 1 for the $s_r$
coordinate of each $u^{(r)}$, for the $t_{r_{d+1}}$ coordinate of
$v^{(d+1)}$,\dots, for the $t_{r_{d+a+1}}$ coordinate of $v^{(d+a+1)}$, and
for the $t_r$ coordinate of each $w^{(r)}$, and we substitute zero
for all other $s$ and $t$  coordinates.  After this substitution,  each of
$v^{(d+1)}$,\dots, $v^{(d+a+1)}$ 
will be $\beta$--orthogonal to all substituted 
vectors except  $u\adi$, but our polynomial
has only degree $a$ in $u\adi$, so the value we get is  zero.
\end{proof}

Since the multiplicity in $s$ of a product of polynomials is the  sum of the 
multiplicities of the factors, 
it follows for $\nu\geq 1$ that the multiplicity in $s$ of any decomposable  
$2\mu$--linear invariant with $\mu>\nu$ is strictly
greater than $\nu$.  
On the other hand, the  multiplicity of $g$ in $s$ is exactly 
$\nu$, since the $2\mu$--linear
monomial 
$s_{1}^{(1)}\cdots s_{\nu}^{(\nu)}
t_{1}^{(\nu+1)}\cdots t_{\nu}^{(2\nu)}
t_{1}^{(2\nu +1)}\cdots t_{1}^{(2\mu)}$
occurs with coefficient 1. We arrive at

\begin{theorem}\label{fulleven}
Let $m>n\geq 2$ both be even.  Then the $m$--linear $O(n)$--invariant 
$g$  
is indecomposable. 
\end{theorem}

\begin{corollary}
Let $m>n\geq 2$ both be even.  
Then the $m$--linear $O(n,\C)$--invariant $\tilde g$ is 
indecomposable in the ring  $R^{O(n,\C)}(\Z)$.
\end{corollary}

\begin{remark}
Theorem~\ref{odd} for $m>n$ follows from 
Theorem~\ref{fulleven}.  
Indeed, identify $O(n-1)$ with the subgroup of $O(n)$ acting on the 
$x$, $y$ coordinate hyperplane in the standard way and fixing $z$.
Then any  $O(n)$--invariant polynomial may be viewed as  an 
$O(n-1)$--invariant polynomial in just 
the $x$ and $y$  variables (regarding the $z$ variables as constants).
View $g$ or $h$ that way, 
and break it up into its multi-homogeneous components.  
One of these is $g_{(n-1)\times m}$  or $g_{(n-1)\times (m-1)}z^{(m)}$, 
which is an indecomposable $O(n-1)$--invariant by Theorem~\ref{fulleven}. It
follows that $g_{n\times 2\mu}$  and $h_{n\times (2\mu+1)}$ are
not in the subalgebra of $R^{SO(n)}$ generated by the elements of degree less
than $2\mu$.  Since no  
$SO(n)$--invariants of degree 1 exist, indecomposability
follows for $h$ as well as $g$.
\end{remark}

\medskip

We now turn to the construction of an indecomposable $m$--linear special orthogonal 
invariant
 for $m\geq n\geq 4$ both even.   Subtract the sum (\ref{DB}) of Proposition~\ref{DBszorzat} from the Pfaffian in Proposition~\ref{Bszorzat} 
and divide by $2^{\mu-\nu+1}$ to get an $m$--linear 
$SO(n,\mathbb C)$--invariant with integer coefficients. 
Call it $\tilde f=\tilde f_{n\times m}\in R_{n\times m}^{SO(n,\C)}(\Z)$,  
noting that  $\tilde f =\frac 12 \left(\tilde g-\tilde h\right)$.
View $\tilde f$  modulo 2 to get an $m$--linear $SO(n)$--invariant 
$f=f_{n\times m}\in R^{SO(n)}_{n\times m}$ with coefficients in $\mathbb F_2$
that is the sum of those $m$--linear 
monomials $\rho$ that, when written in the form 
\eqref{monom},  have even 
$m_r$ for all $r$, with the $x_r$ and $y_r$ occurring in an alternating order, 
and either all $m_r$ being strictly positive and $y_r$ coming first for an odd 
number of lower indices $r$, or a unique $m_r$ being zero.
Invariance of $f$ again follows from that of $\tilde f$ by 
\cite[Lemma 3.2]{df}.

\begin{theorem}\label{evenspecial}
Let $m\geq n\geq 4$  both be even.  Then the $m$--linear $SO(n)$--invariant $f$ is 
indecomposable. 
\end{theorem}

\begin{corollary}
Let $m\geq n\geq 4$ both be even. Then the $m$--linear $SO(n,\C)$--invariant 
 $\tilde f$ is indecomposable in the ring  $R^{SO(n,\C)}(\Z)$.
\end{corollary}

\begin{proof}
Substituting $z$ for $x_\nu$ and $y_\nu$ in an $SO(2\nu)$--invariant  
yields an 
$O(2\nu-1)$--invariant --- this follows easily from Witt's Theorem \cite[Theorem 7.4]{T}. 
Degrees are  not increased.  The image of $f=f_{n\times m}$ is $g=g_{(n-1)\times m}$, which is of the same degree and is indecomposable by Theorem~\ref{odd}.
The image of a  non-trivial decomposition of  $f$ would
be a non-trivial decomposition  of  $g$, so $f$ must also be indecomposable.
\end{proof}

\end{document}